\documentclass{conm-p-l}
\usepackage{amsfonts,latexsym,amssymb,amsmath,algorithmic,epsfig,rotating}
\newtheorem{dfn}{Definition}[section]
\newtheorem{rem}[dfn]{Remark}
\newtheorem{prop}[dfn]{Proposition}
\newtheorem{thm}[dfn]{Theorem} 
\newtheorem{lemma}[dfn]{Lemma}

\newtheorem{ex}[dfn]{Example}
\newtheorem{algor}[dfn]{Algorithm}
\newtheorem{conj}{Conjecture}
 
\newtheorem{nota}{Notation}
 
\renewcommand{\mod}{ \ \mathbf{mod} \ } 
\newtheorem{aph}{Arithmetic Progression Hypothesis (APH)}

\newlength{\poonen}
\newcommand{\thth}{^{\text{\underline{th}}}}

\newcommand{\stst}{{\text{\underline{st}}}}

\newcommand{\np}{{\mathbf{NP}}}

\newcommand{\conp}{{\mathbf{coNP}}}

\newcommand{\feas}{{\mathtt{FEAS}}}

\newcommand{\gln}{\mathbb{G}\mathbb{L}_n}
\newcommand{\glm}{\mathbb{G}\mathbb{L}_m}

\newcommand{\bpp}{{\mathbf{BPP}}}

\newcommand{\crap}{\pp^{\np^\np}}
\newcommand{\pp}{\mathbf{P}}

\newcommand{\am}{\mathbf{AM}}

\newcommand{\pspa}{\mathbf{PSPACE}}
\newcommand{\expt}{\mathbf{EXPTIME}}

\newcommand{\mcyclo}{{\mathtt{TorsionPoint}}}
\newcommand{\hastor}{{\mathtt{HasTorus}}} 
\newcommand{\eps}{\varepsilon}

\newcommand{\Q}{\mathbb{Q}}

\newcommand{\C}{\mathbb{C}}
\newcommand{\N}{\mathbb{N}}
\newcommand{\Z}{\mathbb{Z}}

\newcommand{\Zn}{\Z^n}

\newcommand{\vd}{{\bar{d}}}

\newcommand{\Cn}{\C^n}

\newcommand{\cF}{{\mathcal{F}}}

\newcommand{\Csn}{{(\C^*)}^n}

\renewcommand{\qed}{$\blacksquare$}
\newcommand{\dia}{$\diamond$}

\newcommand{\lcm}{\mathrm{lcm}}

\newcommand{\size}{\mathrm{size}}

\newtheorem{linnik}{Linnik's Theorem}

\begin{document}
\title{Efficiently Detecting Torsion Points and Subtori}  

\author{
J. Maurice Rojas}\thanks{
Department of Mathematics,  
Texas A\&M University
TAMU 3368, 
College Station, Texas 77843-3368,  
USA. 
e-mail: {\tt 
rojas@math.tamu.edu
},  
web page: {\tt 
www.math.tamu.edu/\~{}rojas} \ .
Partially supported by NSF individual grant 
DMS-0211458, NSF CAREER grant DMS-0349309, and Sandia National 
Laboratories.}   

\date{\today}

\dedicatory{To Helaman Ferguson} 

\begin{abstract}
Suppose $X$ is the complex zero set of a finite collection of polynomials in 
$\Z[x_1,...,x_n]$. We show that deciding whether $X$ 
contains a point all of whose coordinates are $d\thth$ roots of unity 
can be done within $\np^\np$ (relative to the sparse encoding), under a 
plausible assumption on primes in arithmetic progression. In particular, 
our hypothesis can still hold even under certain failures of the Generalized 
Riemann Hypothesis, such as the presence of Siegel-Landau zeroes.  
Furthermore, we give a similar {\bf unconditional} complexity upper bound 
for $n\!=\!1$. 
Finally, letting $T$ be any algebraic 
subgroup of $\Csn$ we show that 
deciding $X \;  \text{\raisebox{.5mm}{\scalebox{1}[.5]{$\stackrel{?}
{\supseteq}$}}}\; T$ 
is $\conp$-complete (relative to an even more efficient encoding), 
unconditionally. 
We thus obtain new non-trivial families of multivariate 
polynomial systems where deciding the existence of complex roots can be done 
unconditionally in the {\bf polynomial hierarchy} --- a family of 
complexity classes lying between $\pspa$ and $\pp$, intimately connected 
with the $\pp\!\stackrel{?}{=}\!\np$ Problem. We also discuss 
a connection to Laurent's solution of Chabauty's Conjecture from arithmetic 
geometry.
\end{abstract} 

\maketitle

\section{Introduction} 
While the algorithmic complexity of many fundamental problems in algebraic 
geometry remains unknown, important recent advances have revealed 
that algebraic geometry and algorithmic complexity are closely and subtly 
intertwined. For instance, consider the problem of deciding 
whether a complex algebraic set --- specified as the zero set of a collection 
of multivariate polynomials ---  is empty or not. This is the {\bf complex 
feasibility problem}, $\feas_\C$, and we denote its restriction to 
any family $\cF$ of polynomial systems by $\feas_\C(\cF)$.\\ 
{\bf Note:} The complexity classes we are about to mention are 
reviewed briefly in Section \ref{sec:alg} (see \cite{papa} for 
an excellent introductory account).  

Before seminal work of Pascal Koiran \cite{hnam}, the only connection known 
between $\feas_\C$ and the $\pp\!\stackrel{?}{=}\!\np$ problem 
was that $\feas_\C$ is $\np$-hard, i.e., a polynomial time algorithm  
for $\feas_\C$ would imply $\pp\!=\!\np$. 
(The $\pp\stackrel{?}{=}\np$ problem is the most famous open problem 
from theoretical computer science and has a vast 
literature (see, e.g., \cite{21} and the references in \cite{gj,papa}).)  
However, $\np$-hardness tells us little about what complexity 
class $\feas_\C$ actually belongs to, or how quickly we can anticipate 
solving a given instance of $\feas_\C$. 
Koiran's paper \cite{hnam} was the first to show that the truth of the 
{\bf Generalized Riemann Hypothesis (GRH)} yields the 
implication $\feas_\C\!\not\in\!\pp\Longrightarrow \pp\!\neq\!\np$, and 
\cite{dzh} later showed that this implication could still hold 
even under certain failures of GRH. Furthermore, the underlying 
algorithms are entirely different from the usual techniques of commutative 
algebra (e.g., Gr\"obner bases and resultants) and thus breathe new life into 
an old problem. 

Here we present algorithms revealing new non-trivial families $\cF$ of 
multivariate 
polynomial systems where the implication $\feas_\C(\cF)\!\not\in\!\pp 
\Longrightarrow \pp\!\neq\!\np$ holds {\bf unconditionally}. We also 
present several examples indicating that the algorithms yielding our main 
results may be quite practical. 
In the coming sections, we will detail some of the intricacies behind making 
such algorithms free from unproved number-theoretic hypotheses. 
We begin by stating a number-theoretic hypothesis that is 
demonstrably weaker than GRH. We use $\N$ for the positive integers.   

\begin{aph}
There is an absolute \linebreak 
constant $C\!\geq\!1$ such that for any $x,M\!\in\!\N$
with $x\!\geq\!2^{\log^C M}$,
the set\linebreak $\{1+kM\; | \; k\!\in\!\{1, \ldots,x\}\}$
contains at least $\frac{x}{\log^C(xM)}$ primes.
\end{aph}
Assumptions even stronger than APH are routinely used, and widely
believed, in the cryptology and algorithmic number theory communities
(see, e.g., \cite{miller,miha,koiran,jcs,hallgren}).
In particular, while APH is implied by GRH for the number 
fields $\{\Q(\omega_M)\}_{M\in\N}$, where $\omega_M$
denotes a primitive $M\thth$ root of unity, APH can still hold under certain
failures of the latter hypotheses, e.g., the presence of infinitely 
many zeroes off the critical line \cite{dzh}. 
\begin{thm}
\label{thm:plai}
Suppose $f_1,\ldots,f_k\!\in\!\Z[x_1,\ldots,x_n]$, $x\!:=\!(x_1,\ldots,x_n)$,  
and\linebreak $d_1,\ldots,d_n\!\in\!\N$. Let $\mcyclo$ denote the 
following problem: Decide whether the system of equations 
\[f_1(x)\!=\cdots =\!f_k(x)\!=x^{d_1}_1-1\!=\cdots=\!x^{d_n}_n-1\!=\!0\] 
has a solution in $\Cn$. 
Also let the {\bf input size} of the preceding polynomial system be 
$\left(\sum^k_{i=1} \size(f_i)\right)+\sum^n_{i=1} \size(x^{d_i}-1)$, where 
$\size\left(\sum^m_{i=1}c_i x^{a_{i1}}\cdots x^{a_{in}}_n
\right)\!:=\!\sum^m_{i=1}\log\{(|c_i|+2)(a_{i1}+2)\cdots (a_{in}+2)\}$, 
and let $\mcyclo_1$ denote the restriction of $\mcyclo$ to univariate 
polynomials. Then 
\begin{enumerate} 
\item{$\mcyclo\!\in\!\am$, assuming APH.}  
\item{Unconditionally, $\mcyclo_1\!\in\!\np^\np$ and $\mcyclo_1$ is already 
$\np$-hard.} 
\item{When restricted to {\bf fixed} $n$ {\bf and} 
$d_1,\ldots,d_n$, $\mcyclo\!\in\!\pp$ unconditionally.} 
\end{enumerate} 
In particular, $\mcyclo_1\!\not\in\!\pp \Longleftrightarrow \pp\!\neq\!\np$ 
unconditionally. 
\end{thm} 

\noindent 
Our notion of input size is quite natural: To put it roughly, 
$\size(f)$ measures the amount of ink (or memory) one must use to record the 
monomial term expansion of $f$. Note that the degree of a polynomial 
can be exponential in its input size if the polynomial is sparse, e.g.,  
$\size(11z-2xy^{97}z+x^D)\!=\!\Theta(\log D)$. (We employ the usual 
computer science notations $O(\cdot)$ and $\Omega(\cdot)$ to respectively 
denote upper and lower bounds that are asymptotically true up to a 
multiplicative constant. When both conditions hold, then one writes 
$\Theta(\cdot)$.) Thus, in the miraculous event 
that $\pp\!=\!\np$, our algorithm yielding Assertion (2) above has 
complexity {\bf polynomial} in 
the bit-sizes of the $f_i$ and the {\bf logarithms} of the $d_i$ --- 
a property {\bf not} present in any earlier algorithm for $\mcyclo_1$. 

Alternatively, Theorem \ref{thm:plai} tells us that we can try to prove 
$\pp\!\neq\!\np$ by showing that $\mcyclo_1\!\not\in\!\pp$, thus giving another 
opportunity for algebraic geometry tools for the $\pp\!\stackrel{?}{=}\!\np$ 
problem (see also \cite{mulmuley} for a different approach via 
geometric invariant theory). Indeed, should one eventually prove 
{\bf unconditionally} that $\mcyclo$ lies in the polynomial hierarchy 
then it would be more profitable to proceed with an attempt to prove 
$\mcyclo\!\not\in\!\pp$ rather than $\mcyclo_1\!\not\in\!\pp$ (since 
$\mcyclo$ is at least as hard a problem as $\mcyclo_1$). 
\begin{ex}[A Sparse, but Large, Resultant]
\label{ex:first}
Suppose we would like to know if
$f(x_1):=c_1+c_2x^{a_2}_1+\cdots+c_{m-1}x^{a_{m-1}}_1+c_mx^{D}_1$ 
vanishes at some $M\thth$ root of unity, where $m=\Theta(\log^2 M)$, 
the $c_i$ are integers of absolute value bounded above by $10$, and 
$a_2\!<\cdots<\!a_{m-1}\!<\!D\!<\!M$ are 
positive integers.  The classical resultant for two polynomials in one 
variable (see, e.g., \cite{gkz94})
then tells us that $f$ vanishes at an $M\thth$ root of unity iff the 
determinant of a
highly structured $(D+M)\times (D+M)$ matrix vanishes.
Such a matrix is a special case of what is known as a {\bf quasi-Toeplitz}
matrix.  

The best general algorithms for evaluating such determinants 
yields a randomized bit complexity upper bound of 
$O((D+M)^3 \log^\eta (D+M))$, for some 
absolute constant $\eta\!>\!0$ \cite{emirispan}. (Gr\"obner bases, 
being far more general than what we need, yield a deterministic complexity 
upper bound no better than $(D+M)^{O(1)}$ bit operations
(see, e.g., \cite{lakshman}).) More directly, one could also 
compute the gcd of $f$ and $x^M-1$, but this still leads to a deterministic 
bit complexity upper bound no better than $O(DM)$ (see, e.g., 
\cite[Ch.\ 8]{bpr}). Solving even this special
case of $\mcyclo_1$ within $O((D+M)^\eps)$ bit operations for some 
$\eps\!\in\!(0,1)$ is thus still an open problem. \dia
\end{ex}

While the $\np$-hardness of $\mcyclo_1$ was derived 
earlier by David A.\ Plaisted \cite{plaisted} in a different context, 
our complexity upper bounds 
are new: the best previous upper bounds were $\pspa$ 
\cite{pspace} (unconditionally), $\crap$ \cite{dzh}, or $\am$ \cite{hnam} 
(under successively stronger unproved number-theoretic hypotheses, all 
stronger than APH), following from much more general results. 
It is also interesting to note that $\mcyclo_1$ 
is the same as detecting the vanishing of so-called 
{\bf cyclic resultants}, which arise in dynamical systems and knot theory 
\cite{hillar}. 

Let us now motivate and clarify our use of the term ``torsion point'' by 
showing how our results can also be viewed in the context of 
{\bf Lang's Conjecture} from Diophantine geometry (see, e.g., \cite[Conj.\ 6.3, 
pp.\ 37--38]{lang}).
\begin{nota} 
Throughout this paper, we will let $x^a\!:=\!x^{a_1}_1\cdots x^{a_n}_n$ 
and $m\cdot x\!:=\!(m_1x_1,\ldots,m_nx_n)$, 
where it is understood that $a\!=\!(a_1,\ldots,a_n)\!\in\!\Zn$,\linebreak  
$m\!=\!(m_1,\ldots,m_n)\!\in\!\Csn$, 
and $x\!=\!(x_1,\ldots,x_n)\!\in\!\Csn$.  
Also, given $\vd_1,\ldots,\vd_r\!\in\!\Z^n$, we let 
$T\!\left(\vd_1,\ldots,\vd_r\right)$ 
denote the subgroup of $x\!\in\!\Csn$ satisfying 
$x^{\vd_1}\!=\cdots =\!x^{\vd_r}\!=\!1$. 
We call any point of $\Csn$ with each coordinate a root of unity a 
{\bf torsion point}. 
Finally, for any $g_1,\ldots,g_k\!\in\!\Z[x_1,\dots,x_n]$,
$Z(g_1,\ldots,g_k)$ denotes the zero set of $g_1,\ldots,g_k$ in
$\Cn$. \dia 
\end{nota} 

\noindent
The subgroup $T\!\left(\vd_1,\ldots,\vd_r\right)$ is sometimes known in algebraic
geometry as a \linebreak {\bf subtorus},\footnote{The subtori we consider here 
need {\bf not} be connected.} and the set $m\cdot T\!\left(\vd_1,\ldots,\vd_r\right)$ is 
usually called a {\bf translated} subtorus. The distribution of torsion points 
and subtori on algebraic sets happens to be quite special: a given algebraic 
set will have all its torsion points contained in a subset that is 
a finite union of subtori, each translated by a torsion point. 
This follows from a famous result of 
Laurent \cite{laurent} which was conjectured earlier by Chabauty 
\cite{chabauty}. Explicit bounds on how many torsion points can lie in an 
algebraic set have been given by Ruppert in certain cases \cite{ruppert}, and 
Bombieri and Zannieri in far greater generality \cite{bz}. 

Given these deep results, one may suspect that 
$\feas_\C(\cF)$ can be sped up when the underlying family $\cF$ 
is restricted to problems involving torsions points. Our 
two main theorems show that this is indeed the case. In particular,  
Theorem \ref{thm:torsion} below complements Theorem \ref{thm:plai} by 
examining when an algebraic set contains an entire subgroup worth 
of torsion points, as opposed to a single torsion point. Please note 
that Theorem \ref{thm:torsion} does {\bf not} depend on any unproved 
hypotheses. 
\begin{thm} 
\label{thm:torsion} 
Following the notation above, for any $\ell_1,\ldots,\ell_k\!\in\!\N$, 
$\vd_1,\ldots,\vd_r\!\in\!\Z^n$ and 
$f_{i,j}\!\in\!\Z[x_1,\ldots,x_n]$ with $(i,j)$ ranging over 
$\bigcup^k_{i=1}\{(i,1),
\ldots,(i,\ell_i)\}$, let \linebreak $\hastor$ denote the  problem of 
deciding whether\\
\mbox{}\hfill $T\!\left(\vd_1,\ldots,\vd_r\right)
\!\stackrel{?}{\subseteq}\!Z\!\left(\prod^{\ell_1}_{j=1}f_{1,j},\ldots,
\prod^{\ell_k}_{j=1}f_{k,j}\right)$.\hfill\mbox{}\\  
Then, measuring the  underlying input size {\bf instead} as\\  
\mbox{}\hfill $\left(\sum\limits^r_{i=1} \size\left(\bar{d}_i\right)  
\right)+\sum^k_{i=1}\sum^{\ell_i}_{j=1} \size(f_{i,j})$,\hfill\mbox{}\\ 
we have:
\begin{enumerate} 
\item{$\hastor\!\in\!\conp$, and the restriction of $\hastor$ to 
$n\!=\!1$ is already $\conp$-hard.} 
\item{For {\bf fixed} $n$, $\ell_1,\ldots,\ell_k$, and 
$\vd_1,\ldots,\vd_r$, we have $\hastor\!\in\!\pp$.} 
\end{enumerate} 
In particular, 
$\hastor\!\not\in\!\pp\Longleftrightarrow \pp\!\neq\!\np$. 
\end{thm}

\noindent 
Assertions (1) and (2) of Theorem \ref{thm:torsion}, in the special 
case $n\!=\!1$, were derived earlier respectively in  
\cite{plaisted} and Theorem 2 of the first ArXiV version of \cite{thresh}, 
but with no reference to tori. Note in particular that our first 
notion of size for $\prod^\ell_{j=1}g_j$ can be exponential in 
$\sum^\ell_{j=1}\size(g_j)$ (e.g., take $g_j\!:=\!x_j-1$ for all $j$), 
so Theorem \ref{thm:torsion} uses a much more compact notion of 
input size than Theorem \ref{thm:plai}.

Theorems \ref{thm:plai} and \ref{thm:torsion} can thus be viewed as first 
steps toward an algorithmic counterpart to Laurent's Theorem. In particular, 
having derived nearly tight lower and upper complexity bounds,   
our results allow us to efficiently detect the presence of 
subtori. Determining the actual {\bf exceptional locus} --- i.e., 
the precise finite union of translated subtori containing all the torsion 
points in a given algebraic set --- remains an open problem. 

Laurent's Theorem has since been extended to algebraic groups 
more general than $\Csn$ --- semi-Abelian varieties --- by 
McQuillan \cite{mac}, thus solving the aforementioned Lang Conjecture 
\cite[Conj.\ 6.3, pg.\ 37--38]{lang}. 
For instance, a very special case of McQuillan's more general 
result is the Faltings-Mordell Theorem. A very special case of the 
latter result is the fact that an algebraic curve of 
genus $\geq\!2$, say, defined as the zero set of a bivariate polynomial with 
rational coefficients, has at most finitely many rational points. 

The existence of algorithmic counterparts to these 
more general results is thus a tantalizing 
possibility. An implementable algorithm for finding torsion points on 
Jacobians of algebraic curves of genus $\geq\!2$ has already been detailed 
by Bjorn Poonen \cite{bjorn}, and the complexity appears (but has not yet 
been proved) to be polynomial-time for {\bf fixed} genus \cite{bjorn2}. 
Such a complexity bound, if proved for the sparse encoding, would 
form an intriguing analogue to the polynomiality of $\mcyclo$ 
for {\bf fixed} $n$ and $\bar{d}_1,\ldots,\bar{d}_n$. 

In closing this introduction, let us point out that our improved 
complexity bounds appear to hinge on the highly refined structure of the  
Galois groups underlying our equations: cyclic. In particular, whereas complex 
feasibility for an input system $F$ is (conjecturally) solvable by checking 
the density of primes $p$ for which the mod $p$ reduction of $F$ has a root 
mod $p$ \cite{hnam}, our algorithms instead use a {\bf single} 
well-chosen $p$.  
It is therefore appropriate to formulate the following conjecture, 
based on an observation of Rachel Pries \cite{pries}: 
\begin{conj} 
Suppose $\cF$ is the family of polynomial systems $F$ such that 
$Z(F)$ is finite and the Galois group of $F$ over $\Q$ is dihedral 
or bicyclic. Then the restriction of $\feas_\C$ to $\cF$ lies in 
$\pp^{\np^\np}$ unconditionally. 
\end{conj} 

While the algorithm underlying the general case of Theorem \ref{thm:plai} is 
simpler than that of Theorem \ref{thm:torsion}, the key ideas flow more 
clearly if we begin with the latter theorem. So we review 
some key ideas in one variable in Section \ref{sec:ex}, and then prove 
Theorem \ref{thm:torsion} in Section \ref{sec:alg} below. We finally prove 
Theorem \ref{thm:plai} across Sections \ref{sec:new} and \ref{sec:comp}, 
and briefly discuss some limits to possible improvements in Section 
\ref{sec:final}. 

\subsection{Comparison to Related Results} 
As mentioned before, our main results improve upon 
Koiran's earlier algorithms for $\feas_\C$ \cite{hnam} by 
relaxing, or removing entirely, his assumption of GRH for certain input 
families. Our success in the setting of 
torsion points and subtori can hopefully be extended to 
situations where the underlying Galois groups are more complicated, 
and membership in the polynomial hierarchy was possible only under 
stronger assumptions \cite{hnam,dzh}. We also point out that the work of 
David Alan Plaisted \cite{plaisted} --- which focussed on 
polynomials in one variable --- was a central inspiration behind this 
paper. Our results extend \cite{plaisted} to multivariate polynomials and  
subtori, and suggest the broader context of computational 
arithmetic geometry \cite{jcs}. 

One should also remember earlier work of Grigoriev, Karpinski, and Odlyzko 
\cite{gko}, where it was shown that one can decide if one 
sparse univariate polynomial divides another, within $\conp$, assuming 
GRH. Our Theorem \ref{thm:torsion} can be viewed as an unconditional extension  
of their result to certain multivariate binomial ideals. 
Needless to say, the 
results of \cite{hnam,dzh} contain those of \cite{gko} as special cases, 
but the more general results still depend on unproved number-theoretic 
hypotheses. 

Finally, we point out that as this paper was being completed, the author 
found the paper \cite{filas} during a {\tt MathSciNet} search. In this  
paper, the authors present a {\bf polynomial-time} algorithm (found by 
their referee  
\cite[Thm.\ 3 and Algor.\ A, pp.\ 959--962; Acknowledgements]{filas})   
for deciding whether a sparse univariate polynomial of degree $D$ is divisible 
by the $d\thth$ cyclotomic polynomial for an input integer $d$ {\bf whose 
factorization is known}. (David A.\ Plaisted 
claimed such a result 20 years before \cite{filas}, but without a proof 
\cite[Top of page 132]{plaisted}.) As a consequence, they prove that for a 
{\bf fixed number of monomial terms}, one can restrict to $d$ with prime 
factors bounded above by a constant, and thus one obtains a bona fide 
polynomial time algorithm since such integers can be factored in 
polynomial time.  An analogous speed-up for the restriction of $\mcyclo_1$ 
to a fixed number of monomial terms appears to remain unknown.  

The techniques of \cite{filas} are quite similar to those of 
\cite{plaisted}, with two exceptions: (1) \cite{filas} makes no use of 
certificates in finite fields and (2) \cite{filas} makes clever use of a 
result of Conway and Jones \cite{cj} stating in essence that polynomials 
vanishing at a primitive $d\thth$ root of unity can not be ``too sparse''  
as a function of $d$. 

Our techniques complement the results of \cite{filas} by 
showing that their main problem lies in the polynomial hierarchy 
unconditionally, even when the number of monomial terms varies and 
the factorization of $d$ is unknown. This 
follows directly from our proof of Theorem 
\ref{thm:torsion}, which also extends their context to subtori in higher 
dimensions. 

\section{Roots of Unity, Primes, and Illustrative Examples} 
\label{sec:ex}
\begin{dfn} 
For any ring $R$ we will let $R^*$ denote the group of 
multiplicatively invertible elements of $R$. Also, 
a {\bf primitive} $M\thth$ root of unity is a complex number $\omega\!\in\!\C$ 
such that $\omega^M\!=\!1$ and [$\omega^{M'}\!=\!1 \Longrightarrow M|M'$]. 
The $M\thth$ {\bf cyclotomic polynomial}, $\Phi_M\!\in\!\Z[x_1]$, is then 
the minimal polynomial for the primitive $M\thth$ roots of unity. \dia 
\end{dfn} 

\begin{ex}
\label{ex:91} 
Specializing Example \ref{ex:first} from the Introduction, note that 
the following assertions are equivalent: (1) $f$ 
vanishes at an  $M\thth$ root of unity, (2) $f$ vanishes at a 
{\bf primitive} $d\thth$ root of unity for some $d|M$,  
(3) $\Phi_d(x_1)|f(x_1)$ for some $d|M$. 
For the sake of illustration, let us assume $91|M$ and take $d\!=\!91$. Since 
$x^M_1-1\!=\!\prod_{d|M} \Phi_d(x_1)$ for all $M$ (see, e.g., 
\cite[Ch.\ 6]{bs}), it is 
then  easy to see that $f$ vanishes at a primitive $91^\stst$ root of 
unity $\Longleftrightarrow (x^{91}_1-1)|f(x_1)(x^{13}_1-1)(x^7_1-1)$. 
The latter condition is in turn equivalent to the truth of \\ 
($\star$)\hfill$\left(x^{91c}_1-1\right)|
f\left(x^c_1\right)\left(x^{13c}_1-1\right)\left(x^{7c}_1-1\right)$
\hfill\mbox{}\\ for all $c\!\in\!\N$. \dia
\end{ex} 

Our main algorithmic tricks --- when specialized to the example above  
--- are (a) reducing the last check over all $c\!\in\!\N$ to a 
single well chosen $c$ and (b) working over a finite field instead of 
$\Z[x_1]$. In particular, assuming $91c+1$ is prime, it follows easily 
from Fermat's Little Theorem that $(\star) \ \Longrightarrow$
$f(x^c_1)(x^{13c}_1-1)(x^{7c}_1-1)\!\equiv\!0 \mod 91c+1$ for all 
$x_1\!\in\!(\Z/(91c+1)\Z)^*$. The multivariate lemma below will later help us 
derive that the {\bf converse} holds as well, provided $c$ is large enough. 
\begin{lemma} 
\label{lemma:div} 
For any polynomials $g,g_1,\ldots,g_k\!\in\!\Z[x_1,\ldots,x_n]$ 
(expressed as sums of monomial terms), let $\|g\|_1$ denote the sum of the 
absolute values of the coefficients of $g$, and let $d_i\!:=\!\deg_{x_i} g$ for 
all $i$. Then $\left\|\prod^k_{j=1} g_j
\right\|_1\!\leq\!\prod^k_{j=1}\|g_j\|_1$. Also, if $q$ is a prime satisfying 
$q\!>\!\|g\|_1,1+\max_i\{d_i\}$; and $g(x)\!\equiv\!0 \mod q$ for all 
$x\!\in\!((\Z/q\Z)^*)^n$, then $g$ is identically $0$. 
\end{lemma} 
\begin{rem} 
One should recall {\bf Schwartz's Lemma} \cite{schwartz}, which asserts 
that for any field $K$, and any finite subset $S\!\subseteq\!K$, a 
polynomial $g\!\in\!K[x_1,\ldots,x_n]$ that is not identically 
zero vanishes at $\leq\!(d_1+\cdots+d_n)\#S^{n-1}$ points of 
$S^n$. Applying this result would, however, yield a weaker version of 
the second part of our lemma by requiring a larger $q$ ($q\!>\!\sum_i d_i$). 
Nevertheless, the proof below is quite reminiscent of the proof of  
Schwartz's Lemma. \dia 
\end{rem} 

\noindent 
{\bf Proof of Lemma \ref{lemma:div}:} Writing $g_j(x)\!=\!\sum
\limits_{a\in A_j} c_{j,a} x^a$ for all $j$, observe that\\ 
\mbox{}\hfill $\displaystyle{\left\|\prod^k_{j=1}g_j\right\|_1 = 
\left\| \prod^k_{j=1} \sum_{a\in A_j} c_{j,a}x^a 
\right\|_1 =\sum_{\substack{a=a_1+\cdots+a_k\\a_j\in A_j \text{ for all } j}}
\left | \sum_{\substack{(a'_1,\ldots,a'_k)\in (A_1,\ldots,A_k)\\ 
a'_1+\cdots+a'_k=a}} \prod^k_{j=1} c_{j,a'_j}\right|}$\hfill\mbox{}\\
\mbox{}\hfill
$\displaystyle{\leq \sum_{\substack{a=a_1+\cdots+a_k\\a_j\in A_j 
\text{ for all } j}} \sum_{\substack{(a'_1,\ldots,a'_k)\in (A_1,\ldots,A_k)\\
a'_1+\cdots+a'_k=a}} \prod^k_{j=1} |c_{j,a'_j}|\\
\leq \prod^k_{j=1} \sum_{a_j\in A_j} |c_{j,a_j}|=\prod^k_{j=1} \|g_j\|_1}$.
\hfill\mbox{}\\
So the first portion is proved. 

We now proceed by induction on $n$: 
If $n\!=\!1$ then 
$g(x_1)\!\equiv\!0 \mod q$ for all $x_1\!\in\!(\Z/q\Z)^* 
\Longrightarrow c_0\!\equiv\cdots\equiv\!c_{d_1}\!\equiv\!0 \mod q$, 
since $q-1\!>\!d_1$ and a (not identically zero) polynomial of 
degree $\leq\!d_1$ can have at most $d_1$ roots in $(\Z/q\Z)^*$.
Since $q\!>\!\|h\|_1\!\geq\!\max_i|c_i|$, 
we thus have $c_0\!=\cdots=\!c_{d_1}\!=\!0$, and our base case is 
complete. 

To conclude, assume that the second portion of our lemma holds for 
some fixed $n\!\geq\!1$. Let us then temporarily consider 
$g$ as a polynomial in $x_{n+1}$ with coefficients 
in $\Z[x_1,\ldots,x_n]$. Let $c_i(x_1,\ldots,x_n)$ denote the 
coefficient of $x^i_{n+1}$. Fixing any values for 
$x_1,\ldots,x_n$, observe that just as in the last 
paragraph, $g$ can vanish at no more than $d_{n+1}$ values of 
$x_{n+1}\!\in\!(\Z/q\Z)^*$. Since $q-1\!>\!d_{n+1}$ we then obtain 
$c_0(x_1,\ldots,x_n)\!\equiv\cdots 
\equiv\!c_{d_{n+1}}(x_1,\ldots,x_n)\!\equiv\!0 \mod q$ for all 
$x_1,\ldots,x_n\!\in\!(\Z/q\Z)^*$. Since 
$\|c_i(x_1,\ldots,x_n)\|_1\!\leq\!\|h\|_1$ 
for all $i$, and since the $c_i$ have exponents no larger than 
$q-2$, our induction hypothesis then implies that 
$c_0,\ldots,c_{d_{n+1}}$ are identically $0$, and thus 
$g$ is indeed identically $0$. \qed 

\medskip 

That we can pick a {\bf small} $c$ with $cM+1$ prime is guaranteed by a 
classic theorem of Linnik. 
\begin{linnik} 
The least prime of the form $cM+b$, where $M$ and $b$ are relatively 
prime integers and $1\!\leq\!b\!<\!M$, does not exceed 
$M^{C_0}$ for some absolute constant $C_0$. \qed 
\end{linnik} 

\noindent
The best current unconditional estimate for $C_0$ is $C_0\!\leq\!5.5$, 
assuming  
$M$ is sufficiently large \cite{heathbrown}. 
It is also known that the truth of GRH implies that we can 
take $C_0\!=\!2+\eps$ for any $\eps\!>\!0$,  
but of course valid only for $M\!>\!M_0$, with $M_0$ an increasing 
function of $\frac{1}{\eps}$ \cite[Thm.\ 8.5.8, pg.\ 223]{bs}.  
\begin{ex} 
[A Number-Theoretic Speed-Up] 
\label{ex:biguni}
Let us consider

\noindent
{\scriptsize
\begin{eqnarray*}
f(x_1)& := &x^{255255}_1-5x^{249255}_1-3x^{248928}_1+4x^{234655}_1
-5x^{221135}_1 +2x^{213883}_1-x^{210952}_1+ 4x^{200774}_1\\
& & +4x^{199666}_1 -5x^{191411}_1 +5x^{187436}_1 +2x^{186678}_1 
+3x^{181717}_1-4x^{181453}_1 +5x^{180273}_1+3x^{176054}_1\\
& & +3x^{171282}_1 -4x^{170662}_1+3x^{168177}_1 +x^{164270}_1 
+5x^{157315}_1+2x^{154380}_1 +5x^{147177}_1 -2x^{144498}_1\\
& & -4x^{142969}_1-2x^{139399}_1+3x^{127018}_1+3x^{103857}_1 
-4x^{101698}_1+x^{97641}_1+2x^{91638}_1-5x^{88391}_1\\
& & -5x^{88198}_1-4x^{86818}_1+5x^{85759}_1+5x^{73803}_1-x^{64076}_1
-3x^{60689}_1 -2x^{50793}_1-5x^{24214}_1+4x^{22380}_1\\
& & -2x^{12176}_1-5x^{682}_1-2, 
\end{eqnarray*} }

\noindent
which has degree $255255$ and exactly 
$46$ monomial terms, and 
suppose we'd like to verify whether $f$ vanishes at some 
$510510\thth$ root of unity. To illustrate our approach via 
cyclotomic polynomials, let us first see if $f$ vanishes 
at a {\bf primitive} $91^\stst$ root of unity.  
As observed earlier, when $q\!:=\!91c+1$ is prime, we have that 
$(x^{91c}_1-1)|f(x^c_1)(x^{13c}-1)(x^{7c}-1) \Longrightarrow 
f(t^c)(t^{13c}-1)(t^{7c}-1)\!\equiv\!0 \mod q$ for all 
$t\!\in\!(\Z/q\Z)^*$. So Condition ($\star$) implies a certain 
congruence holds. However, the reduction goes the other way as well: Lemma 
\ref{lemma:div} (applied to the {\bf mod $\pmb{t^{91}-1}$ reduction} of 
$f(t)(t^{13}-1)(t^7-1)$) tells us that the {\bf converse} to the 
preceding implication holds, provided $q$ is prime,  
$q\!>\!\|f\|_1\|x^{13}_1-1\|_1\|x^7_1-1\|_1\!=\!568$, and $q\!>\!255256$.  

In particular, $2842\cdot91 +1\!=\!258623$ is prime. So to check whether 
$f$ vanishes at a primitive $91^\stst$ root of unity, we need only check 
whether \[ f\!\!\left(t^{2842}\right)\left(t^{2842\cdot 13}-1\right)
\left(t^{2842\cdot 7}-1
\right)\stackrel{?}{\equiv}0 \mod 258623 \text{ for all } 
t\in (\Z/258623\Z)^*.\] 
Since $t\!=\!3$ yields $76177$ for the product polynomial above, 
we thus have certification that $f$ does {\bf not} vanish at any 
primitive $91^\stst$ root of unity. 
Similar calculations for 
{\bf small} choices of $c$ and $t$ then suffice to show that $f$ does 
not vanish at any primitive $d\thth$ root of unity for any other  
$d|510510$ either. (Excluding the easy case $d\!=\!1$ and the case $d\!=\!91$ 
we just did, there are exactly $126$ other such cases.) Thus, we can at last 
certify that $f$ does not vanish at {\bf any} $510510\thth$ root of unity. 
\dia 
\end{ex} 

It is easily checked that the number of bit operations 
for the calculations of Example \ref{ex:biguni} (including the work 
for the additional $126$ cases of $d|510510$) lies in the 
lower hundreds of thousands. 
(This is via standard mod $n$ arithmetic (see, e.g., \cite[Ch.\ 5]{bs}), with 
{\bf no} use of FFT multiplication.) More concretely, the finite field 
certificate check above took but a fraction of a second.\footnote{Using the 
computer algebra system {\tt Maple 9.5}, on {\tt diana}, the author's 4Gb 
dual-Athlon 2 Ghz Fedora Core 4 Linux system.} On the other hand, 
computing the gcd of $x^{510510}-1$ and the $f$ above took 
37 minutes and 38.9 seconds.$^2$ We analyze the underlying asymptotic 
complexity in greater depth in the next section, where we also formalize our 
algorithm for $\hastor$. 

\section{Complexity Issues and Proving Theorem \ref{thm:torsion}: 
Detecting Subtori Unconditionally}  
\label{sec:alg} 
Let us recall the following informal descriptions of 
some famous complexity classes. A completely rigourous and 
detailed description of the classes below can be found in the 
excellent reference \cite{papa}. Our underlying computational model 
is the classical Turing model. For concreteness, it is not 
unrealistic to simply imagine that we are working with a laptop computer, 
equipped with infinite memory, flawless hardware, and a flawless 
operating system: classical theorems from complexity theory allow one 
to define the complexity classes below in a machine-independent manner.  
(We omit these more formal definitions for brevity). In particular, 
we can identify ``time'' or ``work'' with how long our laptop computer 
takes to solve a given problem, and ``input size'' can simply be 
identified with the number of bytes in some corresponding input file.  
\begin{itemize}
\item[$\pp$]{ The family of decision problems which can be done within time
polynomial in the input size.\footnote{Note that the underlying polynomial
depends only on the problem in question (e.g., matrix inversion, shortest path
finding, primality detection) and not the particular instance of the
problem.}}
\item[$\bpp$]{ The family of decision problems admitting randomized algorithms
that terminate in polynomial-time to give an answer which is correct with
probability at least\footnote{It is easily shown that we can
replace $\frac{2}{3}$ by any constant strictly greater than $\frac{1}{2}$
and still obtain the same family of problems.} $\frac{2}{3}$.}
\item[$\np$]{ The family of decision problems where a {\tt ``Yes''} answer can
be {\bf verified} within time polynomial in the input size.}
\item[$\conp$]{ The family of decision problems where a {\tt ``No''} answer
can be {\bf verified} within time polynomial in the input size.}
\item[$\am$]{ The family of decision problems solvable by a $\bpp$ algorithm
which has been augmented with exactly {\bf one} use of an oracle in $\np$. }
\item[$\np^\np$]{ The family of decision problems where a {\tt ``Yes''}
answer can be certified by using an $\np$-oracle a number of times
polynomial in the input size. }
\item[$\pp^{\np^\np}$]{ The family of decision problems solvable within time
polynomial in the input size, with as many calls to an $\np^\np$ oracle as
allowed by the time bound. }
\item[$\pspa$]{ The family of decision problems solvable within time
polynomial in the input size, provided a number of processors exponential in
the input size is allowed. }
\item[$\expt$]{ The family of decision problems solvable within time
exponential in the input size.}
\end{itemize} 
\noindent 
The inclusions\\ 
\mbox{}\hfill$\pp\subseteq\bpp\cup\np\subseteq\am\subseteq\conp^\np\subseteq
\pp^{\np^\np}$\hfill\mbox{}\\
and\\ 
\mbox{}\hfill $\pp\subseteq\np\subseteq\np^\np\subseteq\pp^{\np^\np}
\subseteq\pspa\subseteq\expt$,\hfill\mbox{}\\
are fundamental in complexity theory \cite{papa,lab}, 
and the properness of every explicitly stated inclusion above
turns out to be a major open problem (as of late 2007). For instance, while 
we know that $\pp\!\subsetneqq\!\expt$, the inclusion
$\pp\!\subseteq\!\pspa$ is not even known to be proper. 
The first $6$ complexity classes in the list above lie in a family 
known as the {\bf polynomial hierarchy}. It is known that
$\pp\!=\!\np$ implies that the polynomial hierarchy {\bf collapses},
which in particular yields the equalities $\pp\!=\!\np\!=\!\conp\!=\!\am\!=
\!\np^\np\!=\!\pp^{\np^\np}$ \cite[Thm.\ 17.9]{papa}. This standard fact
will be used later.

The structure of our main algorithms depends on a useful number-theoretic 
lemma stated below. In what follows, $e_i$ denotes the 
$i\thth$ standard basis vector of whatever finite-dimensional 
module we are working in. 
\begin{dfn} 
For any $g\!\in\!\Z[x_1,\ldots,x_n]$, let $\bar{g}\!\in\!\Z[x_1, \ldots,x_n]$ 
denote the polynomial obtained by reducing all exponent vectors in 
the monomial term expansion of $g$ modulo the additive subgroup 
$\langle d_1e_1\ldots,d_r e_r\rangle$ of $\Zn$ and collecting terms. \dia 
\end{dfn}
Note that computing $\bar{g}$ is nothing more
than repeatedly applying the substitution $x^{d_i}_i\!=\!1$ (for all 
monomial terms and $i\!\in\!\{1,\ldots,r\}$), and simplifying, until one 
obtains a polynomial with
degree $<\!d_i$ with respect to $x_i$ for all $i\!\in\!\{1,\ldots,r\}$.
Note also that any coefficient of $\bar{g}$ is a sum of coefficients of $g$.
\begin{prop} 
\label{prop:bits} 
For any $g,g_1,\ldots,g_\ell\!\in\!\Z[x_1,\ldots,x_n]$, let 
$m_j$ denote the number of monomial terms of $g_j$ for all $j$. 
Then $\|\bar{g}\|_1\!\leq\!\|g\|_1$, and the monomial term 
expansion of $\overline{\prod^\ell_{j=1}g_j}$
can be computed within\\ 
\mbox{}\hfill 
$O\!\left(\min\left\{\prod\limits^\ell_{j=1}m_j,\prod\limits^r_{i=1}d_i\right\}
\left(\sum\limits^\ell_{j=1}\size(g_j)
+\sum\limits^r_{i=1}\log d_i\right)^2\right)$ 
\hfill\mbox{}\\
bit operations. 
\end{prop} 

\noindent 
{\bf Proof:} The first portion follows directly from the 
definition of $\|\cdot\|_1$ and $\bar{g}$. 

To prove the second portion, note that computing $\bar{g}_j$ consists 
simply of reducing the coordinates of the exponent vectors modulo integers
of size no larger than $\max_i\{\log d_i\}$, and then summing up 
coefficients of monomial terms. So via basic fast finite field arithmetic
(e.g., \cite[Table 3.1, Pg.\ 43]{bs}), this can be done within
$\displaystyle{O\!\left(\sum^r_{i=r}\max\left\{\size(g_j),\log(d_i)\right\}
\log\max \left\{ \size(g_j), \log(d_i)\right\}\right)}$ bit operations. 

Next, note that to compute
$\overline{\prod^\ell_{j=1}g_j}$, we can use the recurrence
$G_1\!:=\!\bar{g_1}$, $G_{j+1}\:=\!\overline{G_j\bar{g}_{j+1}}$,
and stop at $G_{\ell}$. Defining $\kappa_j$ to be  
the maximum bit-length of any coefficient of $\bar{g}_j$, 
the number of bit operations to 
compute $G_2$ is then easily seen to be 
$O^*\!\left(\min\left\{m_1m_2,
\prod^r_{i=1} d_i\right\}\left(\kappa_1+\kappa_2+\sum^r_{i=1}\log d_i\right)
\right)$. (The $O^*(\cdot)$ 
notation indicates that additional factors polynomial in 
$\log \kappa_j$ and $\log\log d_i$ are omitted.) This 
bound is obtained by first computing $\bar{g}_1\bar{g}_2$ by simply 
multiplying all monomials of $\bar{g}_1$ with all 
monomials of $\bar{g}_1$ (using fast arithmetic 
along the way), collecting terms, and then reducing the 
exponents as in the definition of $\overline{(\cdot)}$. 
Continuing inductively, our complexity bound follows directly, 
keeping in mind that $\left\|\overline{\prod^\ell_{j=1}g_j}\right\|_1
\!\leq\!\left\|\prod^\ell_{j=1}\bar{g}_j\right\|_1\!\leq\!
\prod^\ell_{j=1}\|\bar{g}_j\|_1\!\leq\!\prod^\ell_{j=1}
\|g_j\|_1$. \qed 

\begin{lemma}
\label{lemma:final}
Following the notation above, suppose $g\!\in\!\Z[x_1,\ldots,x_n]$, 
\linebreak  
$d_1,\ldots,d_r \!\in\!\N$, $D\!:=\!2+\max\{\max_{i\in\{1,\ldots,r\}} 
\{d_i\},\max_{i\in\{r+1,\ldots,n\}}\{\deg_{x_i} g\}\}$,\linebreak  
$M\!:=\!\left\lceil \frac{\max\{\|g\|_1,D\}}{\lcm_i\{d_i\}}
\right\rceil \lcm_i\{d_i\}$, and assume $c$ is a positive integer 
such that \linebreak $q\!:=\!cM+1$ is prime. Then\\ 
\mbox{}\hfill $T(d_1e_1,\ldots,d_re_r)\subseteq Z(g) \Longleftrightarrow 
\left\{
\begin{matrix}\text{$g\!\left(t^{cM/d_1}_1,\ldots,t^{cM/d_r}_r,
t_{r+1},\ldots,t_n\right) \ \equiv \ 0 \mod q$}\\
\text{for all $t_1,\ldots,t_n\!\in\!(\Z/q\Z)^*$.}\hspace{3.2cm}\mbox{}\end{matrix}
\right.$\hfill\mbox{} 
\end{lemma}

\noindent
{\bf Proof:} Let $J$ denote the ideal 
$\langle x^{d_1}_1-1, \ldots,x^{d_r}_r-1\rangle\!\subset\!\overline{\Q}[x_1,
\ldots,x_n]$. Observe that the primary decomposition of $J$ is clearly 
$\bigcap\limits_{\zeta^{d_1}_1=\cdots=\zeta^{d_r}_r=1}  
\langle x_1-\zeta_1, \ldots,x_r-\zeta_r\rangle$, and each ideal in the 
preceding intersection is prime. $J$ is thus a radical ideal in 
$\overline{\Q}[x_1,\ldots,x_r]$. Now let $I\!:=\!J\cap\Q[x_1,\ldots,x_n]$. 
Before proving our desired equivalence we will need the fact 
that the ideal $I$ of $\Q[x_1,\ldots,x_n]$ is radical as well.  
So let us conclude this necessary digression as follows: 

Suppose $f^k\!\in\!I$ for some $f\!\in\!\Q[x_1,\ldots,x_n]$ and 
$k\!>\!1$. Since $J$ is radical and $J\supseteq\!I$, we then clearly obtain 
the existence of 
$f_1, \ldots,f_r\!\in\!\overline{\Q}[x_1,\ldots,x_n]$ with\\ 
\mbox{}\hfill$f(x)\!=\!(x^{d_1}_1-1)f_1(x)+\cdots+(x^{d_r}_r-1)f_r(x)$.
\hfill\mbox{}\\ 
Letting $G$ denote the Galois group of the coefficients of the $f_i$ over 
$\Q$, let us define $f'_i\!:=\!\frac{1}{\#G}\sum_{\sigma\in G}
\sigma(f_i)$ for all $i\!\in\!\{1,\ldots,r\}$. Observe then that $f(x)$ 
also equals $\sum^r_{i=1}(x^{d_i}_i-1)f'_i(x)$, and thus lies in $I$ as well 
by Galois invariance. So $I$ is radical. 

Returning to our main proof, we now see that:\\  
(A) $T(d_1e_1,\ldots,d_re_r)\!\subseteq\!Z(g) 
\Longleftrightarrow g\!\in\!I$, and \\ 
(B) $\{x^a\; | \; a\!\in\!\{0,\ldots,d_1-1\}\times\cdots\times \{0,
\ldots,d_r-1\}\}$ is a $\Q$-vector space basis for\linebreak
\mbox{}\hspace{.7cm}$\Q[x_1,\ldots,x_r]/I$.\\ 
In particular, $T(d_1e_1,\ldots,d_r e_r)\!\subseteq\!Z(g)$ iff 
$\bar{g}$ is identically zero. 
So it suffices to prove that $\bar{g}$ is identically zero $\Longleftrightarrow 
g\!\left(t^{cM/d_1}_1,\ldots,t^{cM/d_r}_r,
t_{r+1},\ldots,t_n\right) \ \equiv \ 0 \mod q $
for all $t_1,\ldots,t_n\!\in\!(\Z/q\Z)^*$. 
Let $I_{cM}\!:=\!\langle x^{cM}_1-1, \ldots,x^{cM}_r-1\rangle$. 

\noindent 
{\bf ($\pmb{\Longrightarrow}$):} 
By (B), $\bar{g}$ identically zero $\Longrightarrow g\!\in\!I$, and thus 
$g(x^{cM/d_1}_1, \ldots,x^{cM/d_r}_r,x_{r+1},\ldots,x_n)\!\in\!I_{cM}$. 
Since $q\!:=\!cM+1$ is prime, Fermat's Little Theorem implies 
$t^{cM}-1\!\equiv\!0 \mod q$ for all $t\!\in\!\{1,\ldots,cM\}$, so\\ 
\mbox{}\hfill $g(t^{cM/d_1}_1,
\ldots,t^{cM/d_r}_r,t_{r+1},\ldots,t_n)\!\equiv\!0$ mod $q$,  
for all $t_1,\ldots,t_n\!\in\!(\Z/q\Z)^*$.\hfill\mbox{}\\ 
(Remember that we have defined $M$ so that $d_i|M$ for all 
$i\!\in\!\{1,\ldots,r\}$.) 

\noindent 
{\bf ($\pmb{\Longleftarrow}$):} 
By (B), $g-\bar{g}\!\in\!I$ for any $g\!\in\!\Z[x_1,\ldots,x_n]$. So 
we must then have \\
\mbox{}\hfill$g\!\left(x^{cM/d_1}_1,\ldots,
x^{cM/d_r}_r,x_{r+1},\ldots,x_n\right)-\bar{g}\!\left(x^{cM/d_1}_1,\ldots,
x^{cM/d_r}_r,x_{r+1},\ldots,x_n\right)\!\in\!I_{cM}$.\hfill\mbox{}\\ 
We therefore obtain that $g\left(t^{cM/d_1},\ldots,t^{cM/d_r}_r,t_{r+1},
\ldots,t_n\right)\mod q$ for all $t_1,\ldots,t_n\!\in\!(\Z/q\Z)^*  
\Longrightarrow\bar{g}\!\left(t^{cM/d_1}_1,
\ldots,t^{cM/d_r}_r,t_{r+1},\ldots,t_n\right)\!\equiv\!0 \mod q$
for all $t_1,\ldots,t_n\!\in\!(\Z/q\Z)^*$, via another application of 
Fermat's Little Theorem.

Now note that $\left\|\bar{g}\right\|_1\!\leq\!\left\|g\right\|\!\leq\!
M\!\leq\!q$ and\\
\mbox{}\hfill $\deg_{x_i} \bar{g}\!\left(x^{cM/d_1}_1,
\ldots,x^{cM/d_r}_r,x_{r+1},\ldots,x_n\right)\!
\leq\!(cM/d_i)(d_i-1)\!<\!cM\!=\!q-1$\hfill\mbox{}\\ for 
all $i\!\in\!\{1,\ldots,r\}$. Furthermore, $\deg_{x_i}\bar{g}\!=\!\deg_{x_i}g
\!\leq\!D-2\!\leq\!M\!<q-1$ for all $i\!\in\!\{r+1,\ldots,n\}$. 
So Lemma \ref{lemma:div} immediately implies 
that $\bar{g}$ is identically $0$. \qed

\smallskip 
We now state our first main algorithm. 
\begin{algor}[For problem $\hastor$, with simplified subtori, 
unconditionally]\mbox{}\\
\label{algor:subtorus}{\bf Input:} Polynomials 
$f_{i,j}\!\in\!\Z[x_1,\ldots,x_n]$ with $(i,j)\!\in\!\bigcup^k_{i=1} 
\{(i,1),\ldots,(i,\ell_i)\}$, positive integers $d_1,\ldots,d_r$, 
and a suitable value for the constant $C_0$ from Linnik's Theorem.\\
{\bf Output:} A true declaration of whether\\ 
\mbox{}\hfill $Z\!\left(\prod^{\ell_1}_{j=1}f_{1,j},\ldots,
\prod^{\ell_k}_{j=1}f_{k,j}
\right)\!\supseteq\!Z(x^{d_1}_1-1,\ldots,x^{d_r}_r-1)$.\hfill
\mbox{}   

\noindent
{\bf Description:} 
\begin{enumerate} 
\setcounter{enumi}{-1}
\item{Replace each $f_{i,j}$ by $\bar{f}_{i,j}$ (following the notation 
above).} 
\item{
Let $N\!:=\!\max_i\left\{\prod^{\ell_i}_{j=1}\|f_{i,j}\|_1\right\}$ and 
$M\!:=\!\left\lceil\frac{\max\{N,D\}}{\lcm_i\{d_i\}}\right\rceil\lcm_i\{d_i\}$, 
where $D$ is\\ $2+\max\left\{\max_{j\in\{1,\ldots,r\}}
\{d_j\},\max_{(i,j)\in\{1,\ldots,k\}\times \{r+1,\ldots,n\}}
\left\{\sum^{\ell_i}_{s=1}\deg_{x_j} f_{i,s}\right\}\right\}$. }   
\item{Nondeterministically, decide whether there is a $c\!\in\!\N$ with 
$c\!\leq\!M^{C_0}$ and $q\!:=\!cM+1$ prime, a $t\!=\!(t_1,
\ldots,t_n)\!\in\!((\Z/q\Z)^*)^n$, and an $i\!\in\!\{1,\ldots,k\}$, such that 
\[ (\heartsuit_i) \hspace{1in} \prod^{\ell_i}_{j=1}
f_{i,j}\!\!\left(t^{cM/d_1}_1,\ldots,
t^{cM/d_r}_r,t_{r+1},\ldots,t_n\right) 
\not\equiv 0 \ \mod \ q. \hspace{1in} \mbox{}\] } 
\item{If the desired $(c,t,i)$ from Step 2 
exists then stop and output {\tt ``NO''}. Otherwise, stop and 
output {\tt ``YES''}. \dia } 
\end{enumerate} 
\end{algor} 

\noindent 
The adverb ``nondeterministically'' can be interpreted in two ways: 
the simplest is to just ignore the word and employ brute-force search. 
This leads to an algorithm which is dramatically simpler and easier 
to implement than resultants or Gr\"obner bases. All of our examples 
were handled this simple way, and the respective timings were already 
competitive with the latter techniques (cf.\ Examples \ref{ex:first},
\ref{ex:91}, \ref{ex:biguni}, and \ref{ex:multi}).

Alternatively, one can observe that Step 2 is equivalent to   
deciding the truth of a quantified
Boolean sentence of the form $\forall y_1 \cdots \forall y_\nu 
B(y_1,\ldots,y_\nu)$, with $B(y_1,\ldots,y_\nu)$ computable in time polynomial 
in the size of our
initial input. This is clarified in our proof of Theorem \ref{thm:torsion} 
below. 

Before starting our proof, we will need a lemma on integral matrices 
to quantify certain monomial changes of variables.  
\begin{dfn}
\label{dfn:hermite}
Let $\Z^{m\times n}$ denote the set of $m\times n$ matrices
with all entries integral, and let $\glm(\Z)$ denote the
set of all matrices in $\Z^{m\times m}$ with determinant $\pm 1$
(the set of {\bf unimodular} matrices).
Recall that any $m\times n$ matrix $[u_{ij}]$ with
$u_{ij}\!=\!0$ for all $i\!>\!j$ is called {\bf upper triangular}.
Then, given any $M\!\in\!\Z^{m\times n}$, we call an 
identity of the form $UM=H$, with $H\!=\![h_{ij}]\!\in\!\Z^{n\times n}$
upper triangular and $U\!\in\!\glm(\Z)$, a {\bf Hermite factorization}
of $M$. Also, if we have the following conditions in addition: 
\begin{enumerate}
\item{$h_{ij}\!\geq\!0$ for all $i,j$.}
\item{\scalebox{.96}[1]{for all $i$, if $j$ is the smallest $j'$ such that
$h_{ij'}\!\neq\!0$ then $h_{ij}\!>\!h_{i'j}$ for all $i'\!\leq\!i$.}}
\end{enumerate}
then we call $H$ {\bf \underline{the} Hermite normal form} of $M$. \dia
\end{dfn}

A {\bf Smith} factorization is a more refined factorization 
of the form $UMV\!=\!S$ with $U\!\in\!\glm(\Z)$, $V\!\in\!\gln(\Z)$, 
and $S$ diagonal. In particular, if $S\!=\![s_{i,i}]$ and we 
require additionally that $s_{i,i}\!\geq\!0$ and $s_{i,i}|s_{i+1,i+1}$ 
for all $i\!\in\!\{1,\ldots,\min\{m,n\}\}$ (setting $s_{\min\{m,n\}+1,
\min\{m,n\}+1}\!:=\!0$), then such a factorization for $M$ is unique 
and is called {\bf the} Smith factorization. 
\begin{lemma}
\label{lemma:unimod}
\cite{unimod,storjo} 
For any $A\!=\![a_{ij}]\!\in\!\Z^{n\times n}$, the Hermite and Smith  
factorizations of $A$ can be computed within $O(n^4
\log^3(n\max_{i,j}|a_{ij}|))$ bit operations. 
Furthermore, the entries of all matrices in these 
factorizations have bit size $O(n^3\log^2(2n+\max_{i,j}|a_{ij}|))$. \qed 
\end{lemma}

\medskip 
\noindent
{\bf Proof of Theorem \ref{thm:torsion}:} 
Define $X\!:=\!Z\!\left(\prod^{\ell_1}_{j=1}f_{1,j},\ldots,
\prod^{\ell_k}_{j=1}f_{k,j}\right)$ and let us first reduce to the
special case where $\vd_i\!=\!d_ie_i$ for all
$i$: Let $M$ be the matrix whose columns are $\vd_1,\ldots,\vd_r$
and define $x^M\!:=\!(x^{\vd_1},\ldots,x^{\vd_r})$.
An elementary calculation then reveals that if
we have the Smith factorization $UMV\!=\!S\!=:\![s_{i,j}]$ (with 
$S$ having exactly $t$ nonzero entries), then $x^M\!=\!1
\Longleftrightarrow 
(z^{s_{1,1}}_1,\ldots,z^{s_{t,t}}_t)\!=\!(1,\ldots,1)$, 
upon setting $x\!:=\!z^U$. Via Lemma \ref{lemma:unimod}, we see
that this change of variables can be found within $\pp$ and
the increase in our input size is polynomial in $O(\size(\vd_1)+
\cdots+\size(\vd_n))$. So let us assume henceforth that  
$\vd_i\!=\!s_ie_i$ (and let $d_i\!=\!s_i$) for all $i\!\in\!\{1,\ldots,t\}$ 
and set $r\!=\!t$. 

The equivalence of
$\hastor\!\not\in\!\pp$ and $\pp\!\neq\!\np$ follows immediately from our
earlier remarks on the polynomial hierarchy \cite[Thm.\ 17.9]{papa}, assuming 
we indeed have $\hastor\!\in\!\conp$. So let us proceed with proving 
Assertions (1) and (2). 

\noindent
{\bf Assertion (1):} 
The $\conp$-hardness of the $n\!=\!1$ restriction of $\hastor$ --- 
stated equivalently as a problem involving sparse polynomial 
division --- is essentially \cite[Thm.\ 4.1]{plaisted}. 
So we need only show that $\hastor\!\in\!\conp$ for general $n$ and, 
thanks to our preceding reductions, this can be done by proving that 
Algorithm \ref{algor:subtorus} is correct and runs within $\conp$. 

Correctness follows immediately from Lemma \ref{lemma:final} 
applied to the polynomials from 
$\text{(}\heartsuit_1\text{)},\ldots,\text{(}\heartsuit_k\text{)}$. 

To analyze the complexity of Algorithm \ref{algor:subtorus}, first note that 
Steps 0 and 1 can clearly be done in polynomial time and Step 3 takes 
essentially constant time. So it suffices to focus on the 
complexity of Step 2. 

Let us then observe that for any $t_1,\ldots,t_n\!\in\!\Z/q\Z$, we can 
verify ($\heartsuit_i$) in polynomial-time: By basic finite 
field arithmetic (see, e.g., \cite[Ch.\ 5]{bs}), we can clearly decide within 
$\pp$ whether any $f_{i,j}$ vanishes at a given point in 
$(\Z/q\Z)^n$ using a number of bit operations polynomial in 
$\size(g)\log q$, and we then simply multiply the appropriate $f_{i,j}$.  
In total, checking $\text{(}\heartsuit_1\text{)},\ldots,\text{(}\heartsuit_k
\text{)}$ at any given point in $(\Z/q\Z)^n$ requires a number of bit 
operations at worst $k$ times a polynomial in 
\mbox{}\hfill $\log(q)\left[
\left(\sum\limits^r_{i=1} \size(d_i)
\right)+\sum^k_{i=1}\sum^{\ell_i}_{j=1} \size(f_{i,j})\right]$.\hfill\mbox{}\\
Now observe that $\size(q)\!=\!O(\log M)\!=\!O(\log(N) + 
\log(D)+ \sum^r_{i=1} \log d_i)$, which is clearly linear in our input 
size. 
Note also that the integer $N$ from Algorithm \ref{algor:subtorus} 
(which by definition is no larger than $M$) is clearly an upper bound on 
the $1$-norms of the polynomials from ($\heartsuit_1$), $\ldots$, 
($\heartsuit_k$). So any instance of inequality ($\heartsuit_i$) 
can clearly be checked in $\pp$. 

Now note that verifying $q\!=\!cM+1$ is indeed prime 
can be done in time polynomial in $\log q$ (which is in turn 
polynomial in our input size): One can either use the 
succinct primality certificates of Pratt \cite{pratt},  
or the deterministic polynomial-time primality testing algorithm from 
\cite{aks}. So Step 2 is nothing more than verifying the 
truth of the following quantified sentence: \\
\mbox{}\hfill \scalebox{1}[1]{$\exists c \exists  
t_1 \cdots \exists t_n \exists i \left[ 
(cM+1 \text{ prime})\wedge (c\!\leq\!M^{C_0})\wedge 
(\heartsuit_i)\right]$.} \hfill\mbox{}\\ 
$X$ contains the subtorus $T(d_1e_1,\ldots,d_re_r)$ iff the preceding 
sentence is {\bf false}. So via our preceding observations, the truth of 
the sentence being quantified 
can be verified in $\pp$, and our algorithm thus runs in $\conp$. 

\medskip 
\noindent
{\bf Assertion (2):} Suppose $n$, $\ell_1,\ldots,\ell_k$, $d_1,\ldots,d_n$ 
are fixed. Then, by Proposition \ref{prop:bits} (with $\ell$ constant), 
we can decide $\hastor$ in $\pp$ simply by reducing the exponents modulo 
suitable integers and doing a brute-force check of the congruence condition 
given by Lemma \ref{lemma:final}. \qed 

\begin{ex}
\label{ex:trouble} 
While it is tempting to propose a variant of Algorithm \ref{algor:subtorus} to 
detect {\bf translated} subtori, here is an example showing that 
at least one naive extension breaks down quickly: 
Suppose $q\!=\!kd+1$ is prime (we can take $k\!\geq\!2$ and 
arbitrary large by Linnik's Theorem), $g(x)\!:=\!\gamma x^{dq}-1$ 
with $\gamma\!\equiv\!2^d \mod q$ and $\gamma\!\in\!\{1,\ldots,q-1\}$; 
and we want to see if $g$ vanishes at {\bf half} of 
every $d\thth$ root of unity. Since this happens iff $g(x/2)$ vanishes at every 
$d\thth$ root of unity, we could try to mimic Algorithm \ref{algor:subtorus} 
by checking whether 
$g(t^{(q-1)/d}/2)\!\equiv\!0 \mod q$ for all $t\!\in\!(\Z/q\Z)^*$. This is 
indeed so, since $g(t^{(q-1)/d}/2)\!=\!\gamma\cdot 
\frac{t^{(q-1)q}}{2^{dq}}-1\!\equiv\!
2^d\cdot \frac{1}{2^d}-1\!\equiv\!0 \mod q$.  However, 
$g(\zeta/2)\!=\!\frac{\gamma}{2^{(q-1)d}}-1\!\neq\!0$ for $\zeta$ any 
$d\thth$ root of unity. \dia
\end{ex}

\section{From Subtori to Torsion Points: Theorem \ref{thm:plai} in One 
Variable, Unconditionally}  
\label{sec:new} 
With some modifications, Algorithm \ref{algor:subtorus} --- which we used to 
detect subtori --- can be used to efficiently find torsion points in the 
univariate case.  
\begin{algor}[For $\mcyclo_1$, unconditionally]\mbox{}\\
\label{algor:torsion1}
{\bf Input:} Polynomials
$f_1,\ldots,f_k\!\in\!\Z[x_1]$ and a positive integer $d$.\\
{\bf Output:} A true declaration of whether
$Z(f_1,\ldots,f_k)$ contains a point $\zeta$ 
with $\zeta^{d}\!=\!1$. 

\noindent
{\bf Description:}
\begin{enumerate}
\item{Using Algorithm \ref{algor:subtorus}, nondeterministically 
decide whether there is a $\delta|d$ with  
$Z(f_1g_\delta,\ldots,f_kg_\delta)\!\supseteq\!Z\!\left(x^\delta_1-1\right)$ 
where\\
\mbox{}\hfill $g_\delta(x_1)\!:=\!\prod\limits_{p
\text{ a prime dividing } \delta }\left(x^{\delta/p}_1-1\right)$.
\hfill\mbox{}} 
\item{If the desired $\delta$ from Step 1 exists then stop and 
output {\tt ``YES''}.\\ 
Otherwise, stop and output {\tt ``NO''}. \dia } 
\end{enumerate} 
\end{algor} 

\noindent 
Just as in our last algorithm, the adverb ``nondeterministically'' can be 
interpreted in two ways: first, one can simply employ brute-force search, 
and this strategy is dramatically simpler and easier
to implement than resultants or Gr\"obner bases. All of our examples 
were handled in this simple way, and the respective timings were already 
competitive with the latter techniques (cf.\ Examples \ref{ex:first}, 
\ref{ex:91}, \ref{ex:biguni}, and \ref{ex:multi}). 
 
Alternatively, one can observe that Step 1 is equivalent 
to deciding the truth of a quantified
Boolean sentence of the form $\exists y_1 \cdots \exists y_{\nu'} 
\forall y_{\nu'+1} \cdots \forall y_\nu B(y_1,\ldots,y_\nu)$, with 
$B(y_1,\ldots,y_\nu)$ computable in time polynomial in the size of our
initial input. This type of sentence forms one of the definitions 
of the complexity class $\np^\np$. 

\medskip 
\noindent
{\bf Proof of Assertion (2) of Theorem \ref{thm:plai}:} The $\np$-hardness 
of $\mcyclo_1$ is already implicit in the proof of \cite[Thm.\ 5.1]{plaisted}, 
so we need only show that $\mcyclo_1\!\in\!\np^\np$. To do the latter, we  
will prove the correctness of Algorithm \ref{algor:torsion1} and 
that it indeed runs within $\np^\np$.

The correctness of Algorithm \ref{algor:torsion1} follows immediately 
from Step 1 and the correctness of Algorithm \ref{algor:subtorus}. In 
particular, it is clear that $f_i$ vanishes at a primitive 
$\delta\thth$ root of unity (indeed, at {\bf all} primitive 
$\delta\thth$ roots of unity) iff $(x^\delta_1-1)|f(x_1)g_\delta(x_1)$. 

Recalling that we've already proved that Algorithm \ref{algor:subtorus} runs 
in $\conp$ in the last section, Step 1 thus consists of a 
single existential quantifier calling a $\conp$ algorithm. 
In particular, verifying that a putative $\delta$ satisfies $\delta|d$ can 
clearly be done in $\pp$, and thus Algorithm \ref{algor:torsion1} runs 
in $\np^\np$.  \qed 
\begin{rem} 
One can show that the number of possible $\delta$ dividing $d$ 
in Step (1) of Algorithm \ref{algor:torsion1} is $O(d^{\eps})$ 
(for any $\eps\!>\!0$), $O((\log d)^{\log(2)+\eps})$ for a 
fraction of integers approaching $1$ as $d\!\longrightarrow\!\infty$ 
(for any $\eps\!>\!0$), and $O(\log d)$ on average. 
This follows easily from earlier estimates on the number of 
divisors of an integer (see, e.g., \cite{hardy,nicrob,deleglise} 
and the references therein). Practically speaking, this means 
that the main complexity bottleneck in Algorithm \ref{algor:torsion1} is the 
efficient detection of cyclotomic factors. \dia  
\end{rem} 

Before moving to the higher-dimensional case of $\mcyclo$, let us
point out that the product trick underlying Algorithm \ref{algor:torsion1} 
does not naively extend to $n\!>\!1$. 
\begin{ex} 
\label{ex:ref} 
Since $1+\omega_3+\omega^2_3\!=\!0$ for any primitive third root of 
unity $\omega_3$, we see that $1+x+y$ vanishes at a point with coordinates 
third roots of unity.  Can we derive a (polynomial-time certifiable) criterion 
to detect this, in the spirit of Lemma \ref{lemma:div} or 
Step 1 of Algorithm \ref{algor:torsion1}? 

As an initial attempt, one could first consider the product \\
\mbox{}\hfill$(1+x+y)(x-1)(y-1)$\hfill\mbox{}\\
(based on mimicking the use of $f_ig_\delta$ in Algorithm \ref{algor:torsion1}) 
and see if it lies in the ideal $\langle x^3-1,y^3-1\rangle$. The preceding 
product, unfortunately, fails this criterion. 

On the other hand, the larger product\\ 
\mbox{}\hfill$(1+x+y)(1+x+y^2)(x-1)(y-1)$\hfill\mbox{}\\ 
{\bf does} lie in the ideal $\langle x^3-1,y^3-1\rangle$. 
However, the most obvious extension of the latter product results 
in a certificate which can have exponentially many factors in general. \dia 
\end{ex} 

\noindent 
While the latter idea does not obviously yield an efficient 
higher-dimensional extension of Algorithm \ref{algor:torsion1}, 
it does enable one to prove the correctness of a {\bf different} 
(and efficient) higher-dimensional extension of Algorithm \ref{algor:torsion1}. 
This we now detail. 

\section{Completing the Proof of Theorem \ref{thm:plai}} 
\label{sec:comp} 
Let us first state an important quantitative result, which follows 
directly from the effective arithmetic Nullstellensatz of Krick, Pardo, and 
Sombra \cite{kps}.
\begin{thm} 
\label{thm:neg} 
Suppose $f_1,\ldots,f_k\!\in\!Z[x_1,\ldots,x_n]$, $d_1,\ldots,d_n$ are 
positive integers, $F\!:=\!(f_1,\ldots,f_k)$, 
$E\!:=\!\max\{\max_i\deg f_i,\max_i d_i\}$, and $\sigma(F)$ is one plus 
the maximum of the absolute value of the log of any coefficient of 
any $f_i$. Then 
$F(x)\!=\!x^{d_1}_1-1\!=\cdots=\!x^{d_n}_n-1\!=\!0$ has {\bf no} 
complex roots iff there are polynomials $g_1,\ldots,g_k,
h_1,\ldots,h_n\!\in\!\Z[x_1,\ldots,x_n]$, and a positive integer 
$\alpha$, with\\
\mbox{}\hfill ($\star\star$) \ \ \ 
$g_1(x)f_1(x)+\cdots+g_k(x)f_k(x)+h_1(x)(x^{d_1}_1-1)
+h_n(x)(x^{d_n}_n-1)\!=\!\alpha$ \hfill\mbox{}\\ 
identically, and 
\begin{enumerate} 
\item{$\deg g_i,\deg h_i\!\leq\!2n^2E^{n+1}$}
\item{$\log \alpha\!\leq\!2(n+1)^3E^{n+1}(\sigma(F)+\log(k+n)+14(n+1)E 
\log(E+1))$ \qed } 
\end{enumerate} 
\end{thm} 

Since $\alpha$ has no more than $1+\log\alpha$ prime factors, it 
is clear that the identity ($\star\star$) persists --- with 
a {\bf nonzero} right-hand side --- even after reduction 
modulo a prime, for all but finitely many primes. 
This in turn easily implies that lacking torsion 
points (for fixed degree) is a property that persists as one 
passes from $\C$ to most finite fields, and the number of 
exceptions is no more than one plus the right-hand side of Inequality 
(2) above. The following lemma shows how {\bf possessing} torsion 
points persists as one passes from $\C$ to certain special 
finite fields. 
\begin{lemma}
\label{lemma:pos} 
Following the notation of Theorem \ref{thm:neg}, suppose 
$f_1(x)\!=\cdots=\!f_k(x)\!=\!x^{d_1}_1-1\!=\cdots=\!x^{d_n}_n-1\!=\!0$ 
has a complex root. Then the mod $q$ reduction of the preceding system has a 
root in $(\Z/q\Z)^n$ for any $q$ with $q\!\equiv\!1 \mod 
\lcm\{d_1,\ldots,d_n\}$ and $q$ prime. 
\end{lemma} 

\noindent
{\bf Proof:} Letting $F\!=\!(f_1,\ldots,f_k)$, note that $Z(F)$ has a torsion 
point of the specified type iff $Z(F)$ contains a point $\zeta\!=\!(\zeta_1,
\ldots,\zeta_n)$ with $\zeta_i$ a primitive $\delta_i\thth$ 
root of unity, for some positive integers $\delta_1,\ldots,\delta_n$ 
with $\delta_i|d_i$ for all $i$. Note then that the polynomial\\
\mbox{}\hfill $h_i(x_1,\ldots,x_n):=\prod\limits_{\substack{(j_2,\ldots,j_n)\\
j_s \text{ coprime to } \delta_s \forall s\in\{2,\ldots,n\}}} 
f_i(x_1,x^{j_2}_2,
\ldots,x^{j_n}_n)$\hfill\mbox{}\\ 
must satisfy $Z(h_i(x_1,\ldots,x_n)g_{\delta_1}(x_1)\cdots 
g_{\delta_n}(x_n))\!\supseteq\!T(\delta_1e_1,\ldots,\delta_n e_n)$ for all $i$, 
where $g_\delta$ is the polynomial defined in Step 1 of Algorithm 
\ref{algor:torsion1}. 

Now suppose $q\!:=\!c\cdot \lcm\{d_1,\ldots,d_n\}+1$ is prime. 
Then, via the ($\Longrightarrow$) portion of Lemma \ref{lemma:final} 
(which, visible from its proof, does {\bf not} require any 
assumptions on the coefficient size), we 
must have $h_i(x^c_1,\ldots,x^c_n)g_{\delta_1}(x^c_1)\cdots
g_{\delta_n}(x^c_n)$ identically zero on $((\Z/q\Z)^*)^n$. 

Since the roots of $g_{\delta_1}(x^c_1)\cdots
g_{\delta_n}(x^c_n)$ are a proper subset of $((\Z/q\Z)^*)^n$, and 
since $\Z/q\Z$ has no zero divisors, we must 
have that for all $i$, some factor of $h_i$ must have a root in 
$((\Z/q\Z)^*)^n$.  So we are done. \qed  

Our final algorithm is actually the simplest of the three algorithms 
of this paper.
\begin{algor}[For $\mcyclo$ in general, assuming APH]\mbox{}\\
\label{algor:torsion}
{\bf Input:} Polynomials
$f_1,\ldots,f_k\!\in\!\Z[x_1,\ldots,x_n]$, positive integers $d_1,\ldots,
d_n$, and a suitable value for the constant $C\!\geq\!1$ from APH.\\
{\bf Output:} A declaration of whether $Z(f_1,\ldots,f_k)$ contains a 
point $\zeta\!=\!(\zeta_1,\ldots,\zeta_n)$ with $\zeta^{d_i}_i\!=\!1$ 
for all $i$, meaningful and correct with probability $>\!\frac{2}{3}$.

\noindent
{\bf Description:}
\begin{enumerate}

\item{Let $E\!:=\!\max\{\max_i\deg f_i,\max_i d_i\}$, 
$M\!:=\!\lcm\{d_1,\ldots,d_n\}$, and let $\sigma(F)$ be one 
plus the maximum of the log of the absolute value of any coefficient of
any $f_i$.}   
\item{Via recursive 
squaring, find the smallest $J$, $K$, and $L$ such that 
\linebreak $L\!>\!1+2(n+1)^3E^{n+1}(\sigma(F)+\log(k+n)+14(n+1)E
\log(E+1))$, $K\!>\!\max\{e^{C^2},2^{\log^C M},36L^2\log^{2C}M\}$ and 
$J\!>\!\log(6)\log^C(KM)$. } 
\item{Pick no more than $J$ random $j\!\in\!\{1,\ldots,K\}$ until 
one either has $q\!:=\!jM+1$ prime, or $J$ such numbers that 
are all composite. In the latter case, stop and output {\tt ``I HAVE 
FAILED. PLEASE FORGIVE ME.''}.} 
\item{Nondeterministically, decide whether the mod $q$ reduction 
of\\
\mbox{}\hfill 
$f_1(x)\!=\cdots=\!f_k(x)\!=\!x^{d_1}_1-1\!=\cdots=\!x^{d_n}_n-1\!=\!0$\hfill
\mbox{}\\ has a root in $(Z/q\Z)^n$. }  
\item{If there is such a solution then stop and
output {\tt ``YES''}. Otherwise, stop and output {\tt ``NO''}. \dia }
\end{enumerate}
\end{algor}

\noindent 
We are now ready to conclude the proof of Theorem \ref{thm:plai}. 

\medskip 
\noindent
{\bf Conclusion of Proof of Theorem \ref{thm:plai}:} 
The equivalence of
$\mcyclo_1\!\not\in\!\pp$ and $\pp\!\neq\!\np$ follows immediately from our
earlier remarks on the polynomial hierarchy \cite[Thm.\ 17.9]{papa}, assuming
we indeed have $\mcyclo_1\!\in\!\np^\np$. The latter is contained in 
Assertion (2), which we already proved in the last section. 
So let us proceed with proving Assertions (1) and (3). 

\medskip 
\noindent 
{\bf Proof of Assertion (1):} It clearly suffices to show that 
Algorithm \ref{algor:torsion} is correct and runs within $\am$. 

Let $F\!:=\!(f_1,\ldots,f_k)$. 
Correctness follows easily from Theorem \ref{thm:neg} and Lemma 
\ref{lemma:pos}. In particular, observe that $K$ --- the 
size of our sample space of numbers congruent to $1$ mod $M$ --- is 
just large enough so that APH implies $\{1,\ldots,K\}$ contains 
at least $6L$ primes. (This follows easily from the basic 
implication\linebreak  
$x\!\geq\!e^{C^2}\Longrightarrow  \frac{x}{\log^C x}\!\geq\!\sqrt{x}$.) 
Notice also that if $F$ does {\bf not} vanish at any torsion point of 
interest, then the mod $q$ reduction of $F$ {\bf does} vanish at a torsion 
point of interest for at most $L$ primes $q$, thanks to Theorem \ref{thm:neg}. 
So the probability of a false {\tt
``YES'' } answer is $<\!\frac{1}{6}$. Furthermore, 
by a routine binomial probability estimate, using the inequality 
$1-t\!\leq\!e^{-t}$ for $t\!\in\!(0,1)$, we obtain that the 
probability of drawing $J$ composite integers is $<\!\frac{1}{6}$. 
In other words, with probability 
$>\!\frac{2}{3}$, Algorithm \ref{algor:torsion} gives the 
right answer.

To conclude, we need only observe that the seemingly large constants 
nevertheless yield low complexity. In particular, observe that 
the number of random bits necessary to 
do our random sampling is $O(J\log K)\!=\!O\!\left(
[\log^C(M)+\log(L)]^{C+1}\right)$ and the number of 
bit operations we must do is near-linear in $O(J\log K)$ (via 
fast finite field arithmetic \cite[Ch.\ 5]{bs}). 
It is then easily checked that $\log M$ and $\log L$ (and thus $J$) 
are polynomial in our 
input size, so our algorithm is nothing more than a $\bpp$ algorithm, 
followed by a {\bf single} call to an $\np$-oracle. This is 
exactly the definition of an $\am$ algorithm \cite{lab}, so we are done. 

\medskip 
\noindent
{\bf Proof of Assertion (3):} 
Let us fix $n$ and $d_1,\ldots,d_r$, and recall the 
notation of Algorithm \ref{algor:torsion1} and the proof of 
Lemma \ref{lemma:pos}. As observed in the proof of Lemma \ref{lemma:pos}, 
$F$ has a torsion point as specified iff there are positive integers 
$\delta_1,\ldots,\delta_n$ with $\delta_j|d_j$ for all $j$, such that 
for all $i$, the complex zero set of\\ 
\mbox{}\hfill $\left(\prod\limits_{\substack{(j_2,\ldots,j_n)\\
j_s \text{ co-prime to } \delta_s \forall s\in\{2,\ldots,n\}}} 
f_i(x_1,x^{j_2}_2, \ldots,x^{j_n}_n)\right)g_{\delta_1}(x_1)\cdots
g_{\delta_n}(x_n)$\hfill\mbox{}\\
contains $T(\delta_1e_1,\ldots,\delta_ne_n)$. By Lemma \ref{lemma:div}, 
since the number of factors and possible $n$-tuples $(\delta_1,
\ldots,\delta_n)$ is constant, the preceding check can be done 
in $\pp$. \qed 

\begin{ex} 
\label{ex:multi}  
Consider the bivariate polynomial system $F\!:=\!(f,g)$ where 
$f$ and $g$ are respectively\\ 
{\small 
${x}^{3879876}{y}^{4594590}+{x}^{3879876}{y}^{4339335}+{x}^{3879876}{y}^{4084080}
+{x}^{3879876}{y}^{3828825}+{x}^{2909907}{y}^{4594590}+{x}^{3879876}{y}^{3573570}+
{x}^{2909907}{y}^{4339335}+{x}^{3879876}{y}^{3318315}+{x}^{2909907}{y}^{4084080}+
{x}^{3879876}{y}^{3063060}+{x}^{2909907}{y}^{3828825}+{x}^{3879876}{y}^{2807805}+
{x}^{1939938}{y}^{4594590}+{x}^{2909907}{y}^{3573570}+{x}^{3879876}{y}^{2552550}+
{x}^{1939938}{y}^{4339335}+{x}^{2909907}{y}^{3318315}+{x}^{3879876}{y}^{2297295}+
{x}^{1939938}{y}^{4084080}+{x}^{2909907}{y}^{3063060}+{x}^{3879876}{y}^{2042040}+
{x}^{1939938}{y}^{3828825}+{x}^{2909907}{y}^{2807805}+{x}^{3879876}{y}^{1786785}+
{x}^{969969}{y}^{4594590}+{x}^{1939938}{y}^{3573570}+{x}^{2909907}{y}^{2552550}+{
x}^{3879876}{y}^{1531530}+{x}^{969969}{y}^{4339335}+{x}^{1939938}{y}^{3318315}+{x
}^{2909907}{y}^{2297295}+{x}^{3879876}{y}^{1276275}+{x}^{969969}{y}^{4084080}+{x}
^{1939938}{y}^{3063060}+{x}^{2909907}{y}^{2042040}+{x}^{3879876}{y}^{1021020}+{x}
^{969969}{y}^{3828825}+{x}^{1939938}{y}^{2807805}+{x}^{2909907}{y}^{1786785}+{x}^
{3879876}{y}^{765765}+{y}^{4594590}+{x}^{969969}{y}^{3573570}+{x}^{1939938}{y}^{
2552550}+{x}^{2909907}{y}^{1531530}+{x}^{3879876}{y}^{510510}+{y}^{4339335}+{x}^{
969969}{y}^{3318315}+{x}^{1939938}{y}^{2297295}+{x}^{2909907}{y}^{1276275}+{x}^{
3879876}{y}^{255255}+{y}^{4084080}+$\\
${x}^{969969}{y}^{3063060}+{x}^{1939938}{y}^{
2042040}+{x}^{2909907}{y}^{1021020}+{x}^{3879876}+{y}^{3828825}+{x}^{969969}{y}^{
2807805}+{x}^{1939938}{y}^{1786785}+{x}^{2909907}{y}^{765765}+{y}^{3573570}+{x}^{
969969}{y}^{2552550}+{x}^{1939938}{y}^{1531530}+$\\
${x}^{2909907}{y}^{510510}+{y}^{
3318315}+{x}^{969969}{y}^{2297295}+{x}^{1939938}{y}^{1276275}+{x}^{2909907}{y}^{
255255}+{y}^{3063060}+{x}^{969969}{y}^{2042040}+{x}^{1939938}{y}^{1021020}+{x}^{
2909907}+{y}^{2807805}+{x}^{969969}{y}^{1786785}+{x}^{1939938}{y}^{765765}+{y}^{
2552550}+{x}^{969969}{y}^{1531530}+{x}^{1939938}{y}^{510510}+{y}^{2297295}+{x}^{
969969}{y}^{1276275}+{x}^{1939938}{y}^{255255}+{y}^{2042040}+{x}^{969969}{y}^{
1021020}+{x}^{1939938}+{y}^{1786785}+{x}^{969969}{y}^{765765}+{y}^{1531530}+{x}^{
969969}{y}^{510510}+{y}^{1276275}+{x}^{969969}{y}^{255255}+{y}^{1021020}+{x}^{
969969}+{y}^{765765}+{y}^{510510}-{x}^{285285}+{y}^{255255}+2$}\\
\mbox{}\hfill and\hfill \mbox{}\\
{\small ${x}^{4594590}{y}^{285285}+{x}^{4339335}{y}^{285285}-{x}^{4594590}+{x}^{4084080}{y}^{285285}-{x}^{
4339335}+{x}^{3828825}{y}^{285285}-{x}^{4084080}-25\,{y}^{3879876}+{x}^{3573570}{
y}^{285285}-{x}^{3828825}+{x}^{3318315}{y}^{285285}-{x}^{3573570}+{x}^{3063060}{y
}^{285285}-{x}^{3318315}+{x}^{2807805}{y}^{285285}-{x}^{3063060}-25\,{y}^{2909907
}+{x}^{2552550}{y}^{285285}-{x}^{2807805}+{x}^{2297295}{y}^{285285}-{x}^{2552550}
+{x}^{2042040}{y}^{285285}-{x}^{2297295}+{x}^{1786785}{y}^{285285}-{x}^{2042040}-
25\,{y}^{1939938}+{x}^{1531530}{y}^{285285}-{x}^{1786785}+{x}^{1276275}{y}^{285285
}-{x}^{1531530}+{x}^{1021020}{y}^{285285}-{x}^{1276275}+{x}^{765765}{y}^{285285}-
{x}^{1021020}-25\,{y}^{969969}+{x}^{510510}{y}^{285285}-{x}^{765765}+{x}^{255255}{
y}^{285285}-{x}^{510510}+{y}^{285285}-{x}^{255255}-26$,}\\ 
which respectively have degrees $8474466$ and $4879875$, and numbers  
of monomial terms $96$ and $42$. We would like to determine whether 
$F$ vanishes at a point $(\zeta,\mu)$ where both $\zeta$ and $\mu$ are 
$4849845\thth$ roots of unity. 

Algorithm \ref{algor:torsion} tells us we can do so, with a 
controllably small error probability, by finding  
a random prime $q$ of the form $4849845c+1$ and checking if 
the mod $q$ reduction of $F$ has a root in $((\Z/q\Z)^*)^2$. 
Taking $c\!=\!22$ yields the prime $q\!=\!106696591$, and proceeding with 
this choice 
we see that the pair $(75770298,101629661)$ is just such a root. 
This indicates that $F$ may indeed vanish at a pair of $4849845\thth$ roots 
of unity, and running Algorithm \ref{algor:torsion} $r$ times would allow us 
to decide this with an error probability $<\!\frac{1}{3^r}$ by 
taking the answer that occupies the majority. (This example in 
fact vanishes at all $(\zeta,\mu)$ with $\zeta$ and $\mu$ 
{\bf primitive} $95\thth$ roots of unity so, since $95|4849845$, our 
putative answer is correct.) 

One could instead try to compute a Gr\"obner basis for the 
ideal\linebreak $\langle f,g,x^{4849845}-1,y^{4849845}-1\rangle$. 
The resulting basis will then be $\{1\}$ iff $F$ does {\bf not} 
vanish at any pair of $4849845\thth$ roots of unity. Trying one of 
the best Gr\"obner basis engines ({\tt Singular}, version {\tt 3-0-2}), 
we are immediately thwarted: the maximum allowed exponent size is $65536$. 
Trying three smaller examples with respective total degrees $92114$ and $65296$ 
(and respective numbers of monomial terms $70$ and $40$) resulted 
in {\tt ``Out of memory''} errors within about $14$ minutes in all 
cases. 

On the other hand, while a brute-force implementation of Algorithm 
\ref{algor:torsion} can run slowly, the corresponding {\tt Maple} 
implementation has no memory problems for our examples here. \dia  
\end{ex}

\section{Is $\mcyclo$ $\np$-complete? } 
\label{sec:final} 
We close this paper by observing a possible speed-up to our last 
algorithm: One could instead simply attempt to nondeterministically 
guess a {\bf small} number of suitable primes (instead of randomly sampling 
a large set), and then check 
nondeterministically whether one has torsion points modulo these primes.  
In particular, if the number of such ``guessed'' primes is polynomial 
in the input size, then it can be proved via the techniques of this 
paper that such an approach would yield $\mcyclo\!\in\!\np$. 

However, it is not clear how to  
prove that a small enough number of primes can be used. 
In particular, our final example shows that one definitely 
needs to use at least $3$ primes, already for one variable. 
\begin{ex} 
Taking\\ 
\mbox{}\hfill 
$f(x_1)\!:=\!4-3x^{10}_1-3x^{18}_1+4x^{42}_1-x^{60}_1+6x^{81}_1
+5x^{95}_1-2x^{102}_1+3x^{105}_1$\hfill\mbox{}\\ and $d\!:=\!210$, it can be 
checked via {\tt Maple 9.5} (within 2 hours, 13 minutes, and 42.63 seconds) 
that $f$ does {\bf not} vanish at any $d\thth$ root of unity. 
One would prefer to do this check modulo an intelligently chosen 
prime of the form $210c+1$ instead. However, there are 
exceptional primes which, using this approach, would cause one to falsely 
declare that $f$ {\bf does} vanish at a $d\thth$ root of unity. 
In this case it easily checked that the exceptional primes are 
exactly the divisors of the resultant of $f$ and $x^{210}-1$, 
which (up to sign) is\\{\small 
\mbox{}\hfill\scalebox{.72}[1]{$2227699600874096872564585144832612236369963246002360338615319497424201747782488174224095731882015016718028$}\hfill\mbox{}}\\
and factors as\\
{\small 
\mbox{}\hfill\scalebox{.57}[1]{$
(2)^2(13)(29)(37)(43)(61)(71)(1801)(108557)(659101)(69529066111)
(261727038763)(20353321490154047885351)(4491828078538834477370467060773855421)
$.}}
\hfill\mbox{}\\
In particular, we see that the $11\thth$ and $14\thth$ prime factors 
above are both congruent to $1$ mod $210$, and could thus lead 
to false {\tt ``YES''} answers. \dia 
\end{ex} 

It is interesting to note that Pascal Koiran has already given some 
evidence that it may be hard to prove that the more general problem 
$\feas_\C$ is $\np$-complete. His evidence is based on the fact 
that $\feas_\C$ contains a hard circuit-theoretic 
problem \cite[Sec.\ 6]{hnam}. However, such a reduction does not appear to 
be known for  $\mcyclo$, so there may be more hope that $\mcyclo\!\in\!\np$ 
than $\feas_\C\!\in\!\np$. 

\section*{Acknowledgements} 
This paper began while the author was attending 
the PIMS workshop on resolution of singularities, factorization of birational 
mappings, and toroidal geometry (December 11--16, 2004, Banff, Canada). The 
author thanks the organizers (Dan Abramovich, Ed Bierstone,  
Dale Cutkosky, Kenji Matsuki, Pierre Milman, and Jaroslaw 
Wlodarczyk) for their kind invitation. This paper was then 
presented at MEGA 2005 (Porto Conte, Alghero, Sardinia, Italy, June 1, 2005), 
and  completed at Sandia National Laboratories. So the author thanks 
Patrizia Gianni, Tomas Recio, Carlo Traverso, Mark Danny Rintoul 
III, Philippe Pebay, and David Thompson for their generous invitations and 
wonderful choice of settings.  

Special thanks go to Matt Baker for reading the ArXiV version of my 
paper the day it was posted and then pointing out the reference 
\cite{bz}. Thanks also to Jeff Achter, Eric Bach, Alex Buium (who pointed out 
the reference \cite{bjorn}), Petros Drineas, Sidney W.\ Graham, Bjorn Poonen, 
Rachel Pries, and Igor Shparlinski for some nice conversations and/or e-mail 
exchanges. I also thank an anonymous referee for pointing out Example 
\ref{ex:ref}, thus spotting an error in an earlier version of Theorem 
\ref{thm:plai}. 

Finally, I would like to dedicate this paper to Helaman Ferguson, whom I 
had the honor of meeting at the January 2005 MAA-AMS-ASL joint meeting in 
Atlanta. Helaman: thank you for the $\pi$-disk! 

\bibliographystyle{amsalpha}

\end{document}